\documentclass[english,11pt]{article}
\usepackage[T1]{fontenc}
\usepackage{amsthm}
\usepackage{color}
\usepackage{euscript}
\usepackage{textcomp}
\usepackage{subfig}
\usepackage{amsmath}
\usepackage{enumitem}
\usepackage{amssymb}
\usepackage{amsthm}
\usepackage[english]{babel}
\usepackage{graphicx}
\newtheorem*{definition*}{Definition}
\newtheorem{definition}{Definition}
\newtheorem{theorem}{Theorem}

\newtheorem{proposition}[theorem]{Proposition}

\newtheorem*{remark}{Remark}
\newtheorem*{theorem*}{Theorem}

\def\C{{\mathbb C}}

\def\R{{\mathbb R}}
\def\GL{\mathrm{GL}}
\def\Per{\mathrm{Per}}
\def\Res{\mathrm{Res}}
\def\SL{\mathrm{SL}}
\def\Sp{\mathrm{Sp}}
\def\Z{{\mathbb Z}}
\def\const{\mathrm{const}}
\def\cH{\mathcal{H}}
\def\cL{\mathcal{L}}
\def\cM{\mathcal{M}}

\def\cR{\mathcal{R}}

\author{Igor Krichever\thanks{Skolkovo Institute of Science and Technology, National Research University Higher School of Economics, Columbia
		University, e-mail: krichev@math.columbia.edu},
	Sergei Lando\thanks{National Research University Higher School of Economics, Skolkovo Institute of Science and Technology, e-mail: lando@hse.ru},
	Alexandra Skripchenko\thanks{National Research University Higher School of Economics, Skolkovo Institute of Science and Technology, e-mail: sashaskrip@gmail.com}}

\title{Real-normalized differentials\\
	with a single order~$2$ pole}

\begin{document}

\maketitle

\begin{flushright}
	{\it To the memory\\
	of our dear friend and colleague\\
Boris Dubrovin}
\end{flushright}

\begin{abstract}
	A meromorphic differential on a Riemann surface is said to be
{\it 	real-normalized} if all its periods are real. Real-normalized differentials
	on Riemann surfaces of given genus with prescribed orders of their poles
	form real orbifolds whose topology is closely related to that
	of moduli spaces of Riemann surfaces with marked points. Our goal is to develop tools to
	study this topology. We propose a combinatorial model for the real normalized differentials with a single
order $2$ pole and use it to analyze the corresponding absolute period foliation.
\end{abstract}

\section{Introduction}
The general concept of real-normalized differential was introduced in \cite{kr87} in the framework of the spectral theory of periodic linear operators.
\begin{definition}
A meromorphic differential on a Riemann surface is said to be \emph{real-normalized\/} if all its periods are real.
\end{definition}
By itself, this concept is almost equivalent to the concept of a harmonic function on a punctured Riemann surface. Indeed, by the definition the imaginary part
$$y(p):={\Im} \, \int^p d\zeta$$
of the abelian integral of a real normalized differential $d\zeta$ is a {\it single-valued} harmonic function (defined up to adding a constant).

Conversely, let $y(p)$ be a real-valued harmonic function on a Riemann surface with punctures; then locally there exists a unique up to an additive constant conjugate harmonic function $x(p)$. Hence, $y(p)$ uniquely defines the differential $d\zeta=dx+idy$, which, by construction, is a real-normalized holomorphic differential on the complement to the marked points. One can specify asymptotic behavior of $y(p)$ near a marked point by the requirement that $d\zeta$ extends to a meromorphic differential on the Riemann surface and has a fixed principal part at the marked points. Recall that
\begin{definition}
The {\em principal part} of a meromorphic differential at a point $p$ on a Riemann surface $C$ is an equivalence class of meromorphic differentials $\omega$ in a neighborhood of $p$, with the equivalence $\omega\sim \omega'$ if and only if $\omega'-\omega$ is holomorphic at $p$.
\end{definition}
Non-degeneracy of the imaginary part of the Riemann matrix of $b$-periods of a normalized holomorphic differential on a smooth genus $g$ algebraic curve implies that

\begin{proposition}\label{1} For any fixed principal parts of poles with pure imaginary residues, whose sum is zero, there exists a unique real-normalized meromorphic differential $d\zeta$, having prescribed principal parts at the marked points.
\end{proposition}
The real normalized differentials are central in the Whitham perturbation theory of algebraic-geometrical solutions of the integrable systems \cite{kr-aver,kr-tau}.

In \cite{GrKr09} it was shown that certain structures and constructions of the Whitham theory can be instrumental in understanding the geometry of the moduli spaces of Riemann surfaces with marked
points. In particular, a new proof of the Diaz’ bound on the dimension of complete subvarieties of the moduli spaces was obtained. In \cite{Kr11} the real normalized differentials were used for the proof of Arbarello's conjecture~\cite{A74}: {\it any compact complex cycle in ${\mathcal M}_g$ of dimension $g-n$ must intersect the stratum of smooth genus $g$ algebraic curves having a Weierstrass point of order at most $n$}.

The main goal of this paper is to propose a combinatorial model for the moduli space of compact Riemann surface~$S$ with one marked point~$p\in S$ and a chosen $1$-jet~$z$ of a local coordinate at~$p$  via identification of this space with the moduli space of real-normalized differentials having a single pole of order $2$.

\begin{remark}
The $1$-jet of a local coordinate~$z$ at a point of a Riemann surface~$S$ can be interpreted
as a nonzero cotangent vector at this point. Indeed, this cotangent vector is nothing but the one
given by the $1$-form $dz$.
In contrast, the principal part of a meromoprhic differential~$\phi$ at a pole~$p$
of order~$2$, where the residue of~$\phi$ vanishes, $\Res_p\phi=0$, can be identified with a tangent
vector at this point: the value of this tangent vector on a cotangent vector of the form
$df$, with~$f$ vanishing at~$p$, is defined to be $\Res_pf\phi$.
\end{remark}

Any space of differentials is foliated into loci of differentials
having the same periods (absolute period foliation). In holomorphic
dynamics, the study of this foliation plays a crucial role in
understanding the behavior of the Teichm\"uller flow on Teichm\"uller
spaces. This flow preserves the stratification of the space of differentials
by the orders of their zeroes.

In the Whitham theory, it was revealed (\cite{dub,kr-lax,kr-tau}) that the leaves of the isoperiodic foliation of meromorphic differentials admit the structure of Frobenius manifolds, which is a geometric manifestation of the topological quantum field theory.

Our goal is to study the stratification and the foliation of the space of real-normalized differentials with a single order 2 pole using their combinatorial presentation by means of so called cut diagrams.

The paper is organized as follows. In Sec.~\ref{sec:rel}
we recall main known results concerning the absolute period foliation
of spaces of holomorphic differentials. These results are to be compared
with the ours, and can be used to analyze the similarity and the difference in
the properties of REL foliations in holomorphic and real-normalized cases.
In Sec.~\ref{sec:comb} we introduce a combinatorial model for the principal
stratum in the space of real-normalized differentials with a single order~$2$ pole based on the behavior of separatrices of a line field associated to the differential.
 The principal stratum consists of differentials having simple zeroes and such that
none of the zeroes belongs to a separatrix entering another zero.
In Sec.~\ref{sec:princ} we investigate density and connectedness
properties of the constant period loci in the principal stratum
of the space of real-normalized differentials with a single order~$2$ pole.
Sec.~\ref{sec:str1} is devoted to a combinatorial description of the
stratification of the space of real-normalized differentials with a single order~$2$ pole on elliptic curves.

The authors are grateful to E.~Lanneau, F.~Ygouf, G. Calsamiglia, B. Deroin, L. Arzhakova and M.~Nenasheva (Dudina)
for useful communications, and to V.~Zhukov for improving the illustrations. 
S.~Lando and A.~Skripchenko appreciate the support of RSF--ANR Grant, Project
20-41-09009.

\section{REL foliations: historical notes}\label{sec:rel}

The space of real-normalized differentials with a single pole of order~$2$
on genus~$g$ curves is foliated:
a leaf of this foliation passing through a real-normalized Abelian differential~$\phi$
is the connected component, which contains~$\phi$,
in the locus of all differentials whose subgroup of periods coincides with $\Per(\phi)\subset\R$.
For holomorphic differentials, a similar foliation into leaves associated to fixed subgroups $\Per(\phi)\subset\C$
has been extensively studied. In this section, we present a brief description of the results of these studies.
The spaces of holomorphic differentials are stratified with respect
to the orders of the zeroes of the differentials, and
the corresponding foliations are restricted as well to the strata.
A number of the results about the structure of the foliations are related to this stratification.

Denote by $H(d_1, \cdots, d_m)$, $d_1+\dots+d_m=2g-2$, the subvariety in the space of holomorphic differentials
consisting of differentials having zeroes of orders $d_1,\dots,d_m$. Below, we will also use multiplicative
notation for partitions, so that $1^{k_1}2^{k_2}\dots$
is the partition having $k_1$ parts equal to~$1$, $k_2$ parts
equal to~$2$, and so on.
In particular, the {\it principal stratum} $H(1,\dots,1)=H(1^{2g-2})$ consists of generic differentials, that is,
differentials all whose zeroes are simple.

Each stratum $H(d_1, \cdots, d_m)$ admits a natural action of the group $\SL(2,\mathbb R)$ which generalizes the action of $\SL(2,\mathbb R)$ on the space $\GL(2,\mathbb R)/\SL(2,\mathbb Z)$ of flat tori; this action can be described in terms of the \emph{period map} or \emph{period coordinates}. Let $(C, \phi)$ be a point of the stratum $H(d_1, \cdots, d_m)$, $\Sigma\subset C$ be the set of zeroes of $\phi$, $|\Sigma|=m$ and let $\gamma_1,\cdots,\gamma_k$ be a basis of the relative homology group
$H_1(C,\Sigma)$. Then the periods
\begin{equation}
\left(\int_{\gamma_1}\phi, \cdots, \int_{\gamma_k}\phi\right)
\end{equation}
define a period map $\Phi: H(d_1, \cdots, d_m)\to \C^k,$ which is
a system of local coordinates in a neighborhood of $(C,\phi)$.
The choice of the basis $\gamma_1,\dots,\gamma_k$ defines a
covering of the stratum $H(d_1, \cdots, d_m)$ by the
corresponding Teichm\"uller space $\widehat H(d_1, \cdots, d_m)$.

The local coordinate system given by a period map can be written as a $2\times k$ real matrix; so, the $\SL(2,\R)$-action
in these coordinates is linear. The restriction of this action to some fundamental domain of the mapping class group leads to the notion of \emph{Kontsevich--Zorich cocycle}.

Each stratum of the Teichm\"uller space $\widehat H(d_1, \cdots, d_m)$ carries a natural measure $\hat\lambda$ obtained by pulling back the Lebesgue measure on $H^1(C, \Sigma, \mathbb{C})$ by the period map. This measure descends to the measure $\lambda$ on the stratum $H(d_1, \cdots, d_m)$ of moduli space. However,
the volume of the stratum with respect to this measure is infinite. The natural way to solve this issue is to restrict ourselves to the subset of area 1 surfaces. Namely, the \emph{area map} $\alpha$ is defined as follows:
\begin{equation}
\alpha(C, \phi) = \frac{i}{2}\int_{C}\phi\land\bar\phi.
\end{equation}
Desintegration of the measure $\lambda = \int_{\mathbb R_{+}}\lambda_a$ gives a family of measures $\lambda_a$. The measure $\lambda_1$ defined on $\alpha^{-1}(1)$ is called the \emph{Masur--Veech measure}. This measure is known to be ergodic and finite.

Now, each  stratum $H(d_1, \cdots, d_m)$ carries a natural holomorphic foliation:
a leaf of this foliation passing through an Abelian differential~$\phi$
is the connected component, which contains~$\phi$, in the locus of all differentials
whose subgroup of absolute periods coincides with $\Per(\phi)\subset\C$.
This foliation extends to the closure $\overline{H(d_1, \cdots, d_m)}$ of the stratum (so that the intersection of a leaf
in $H(d_1, \cdots, d_m)$ with a more degenerate stratum
is a disjoint union of leaves of the foliation in the latter). This foliation can be considered as a quotient of the foliation defined on the Teichm\"uller space $\widehat H(d_1, \cdots, d_m)$, whose leaves are known to be invariant under the action of the mapping class group.

This foliation was described many times independently (see \cite{McM07_1}, \cite{McM07_2}, \cite{Zo}, \cite{Sch}) and is known under several different names, among which absolute period foliation, isoperiodic foliation, REL foliation and kernel foliation. This object played an important role in the study of dynamics of the action of $GL(2,\mathbb{R})$ on the moduli space. It is known to be invariant under the action of the Teichm\"uller flow and the holonomy map.
The leaves of the absolute period foliation
are transverse to the orbits of the $\SL_2(\R)$-action above;
the dimension of the latter leaves
is $2g-3$ (for $g \ge 2$).

The most interesting questions about REL foliation concern its ergodic properties (transitivity and ergodicity) as well as geometry of the leaves. By \emph{transitivity} we mean the existence of a dense leaf. A smooth foliation is said to be \emph{ergodic} with respect to a given  measure if every Borel set that is a union of its leaves has either a full measure or measure zero.


The first result concerning ergodicity of REL foliation in principal strata is due to McMullen who proved in \cite{McM14} that for every genus~$g$ any fiber of the period mapping on $\cH_g$ is a slice of the Schottky locus in the Siegel space $\cH_g$ of symmetric $g\times g$ matrices with positive definite imaginary part, by a linear copy of the Siegel space $\cH_{g-1}$ inside $\cH_{g}$ (modulo partially compactifying the fibers by adding nodal abelian differentials of compact type) (Theorem 1.1 in \cite{McM14}). In case of $g\le 3$,  this implies that the absolute period foliation is ergodic and its closed leaves all come from elliptic cohomology classes (Proposition 2.6). Recall that a cohomology class $[\phi]\in H^1(C)$ is called \emph{elliptic of degree $d>0$}, if $\Per(\phi)\cong \Z^2$ is a lattice in $\C$ and if the natural map from~$C$ to the elliptic curve $\C/ \Per(\phi)$ has degree $d$. The proof of the ergodic statement is based on the standard tools of homogeneous dynamics (Moore's ergodic theorem for the first statement and Ratner's theorems for the second one) and uses the fact that in case of genus $2,3$ the Schottky locus is the whole Siegel space. As for the closed leaves, Ratner's theorem implies that the leaf of the foliation $\EuScript F$ for a given $\Sigma_g$ defined by $\phi: H_1(\Sigma_g, \mathbb R)\to\mathbb C$ is closed iff $\Gamma = G\cap Sp_{2g}(\mathbb Z)$ is a lattice in $G = Sp_{2g}(\mathbb R)^{\phi}.$

However, for higher genera $(g\ge 5)$, there is no simple explicit way to understand, for a given $\tau\in\mathcal{H}_g$, whether $[\tau]\in\mathcal{A}_g$ belongs to $\mathcal{J}_g$ or not, where $\mathcal{A}_g$ denotes the moduli stack of principally polarized (complex) abelian varieties of dimension $g$ and $\mathcal{J}_g$ is the image of the moduli stack of smooth (compact complex) curves of genus $g$ under the Jacobian map. As a consequence, McMulllen's approach to ergodic properties of REL foliations is hard to generalize.

However, a much more general result is established in \cite{Ha} and \cite{CaDeFr}:
\begin{theorem}
The absolute period foliation of the principal stratum is ergodic in every genus $g\ge 2$.
\end{theorem}

The two proofs of this theorem are morally significantly different. Calsamiglia, Deroin and Francaviglia in \cite{CaDeFr} classify completely the closures of the leaves of the absolute period foliation of the principal stratum and obtain the statement about ergodicity as a consequence of this classification. More precisely, first they prove the so-called \emph{transfer principle}, which asserts  that the fibers of the period map
are connected, and then apply Ratner's theory to the linear action of the integer symplectic group $\Sp(2g,\mathbb Z)$ on the subset in
$\C^{2g}$ corresponding to the periods of abelian differentials. The idea to apply Ratner's theory in this context is due to M.~Kapovich~\cite{Ka}.

The proof suggested by U. Hamenst\"adt in \cite{Ha} is a direct proof of ergodicity statement performed by induction based on McMullen's result in genus 2.

Almost nothing is known about ergodic properties of the absolute period foliation on irreducible components of non-principal strata. The main difficulty comes from the fact that in other strata the fibers of the period map are not necessarily connected (see \cite{Yg} for the concrete examples).

One of the most advanced known results in this direction is due to Hooper and Weiss (see \cite{HoWe}):
\begin{theorem}\label{HW}
The absolute period foliation in the odd connected component of the stratum $\mathcal{H}(g-1, g-1)$ which contains the so-called Arnoux--Yoccoz surface always has a dense leaf.
\end{theorem}

The proof relies on two observations that rarely can be combined: on one hand, Arnoux--Yoccoz surfaces are periodic under the Teichm\"uller flow (because of the pseudo-Anosov property); this implies, in particular, that the horizontal foliation is uniquely ergodic; at the same time, it is known that an arbitrarily small deformation of these surfaces results in a surface such that the horizontal foliation is periodic (the latter property is related to the fact that a famous invariant, defined for a measured foliation on surfaces, named SAF, vanishes for Arnoux--Yoccoz surface).

Theorem \ref{HW} was recently extended in \cite{Yg} to an infinite family of translation surfaces endowed with a zero SAF pseudo-Anosov flow introduced by Do and Schmidt in \cite{DoSch}:
\begin{theorem}
The leaves of the Do--Schmidt surfaces are dense in the stratum in which they belong.
\end{theorem}
Transitivity very likely implies ergodicity, but there is no rigorous proof of any ergodic statement yet.

Finally, a series of results about the isoperiodic dynamics in affine manifolds of rank 1 (mainly, Prym eigenform loci) was established in \cite{Yg}. In particular, the following statement was proved:
\begin{theorem}
Let $M$ be a rank one affine manifold. Then either all the leaves of
the isoperiodic foliation on $M$ are closed or all the leaves are dense in the iso-area
sublocus of $M$ in which they belong. In the second case, the foliation is ergodic with
respect to any $\SL_{2}(\R)$-invariant ergodic measure supported on $M$.
\end{theorem}

It was used to get information about the isoperiodic dynamics in the strata (the problem we started with) for some particular surfaces with foliations on them (namely, for certain classes of Prym eigenforms).
Again, it was shown that the considered foliations are transitive.
These observations were also used for the classification of affine manifolds of genus $g=3$.

We complete the description of the abelian differential case with a brief outline of the notion of the \emph{Rel flow} which seems to be a very efficient tool to study the absolute period foliations (see \cite{Yg} and \cite{HoWe} for the details). For a given point of moduli space $(C,\phi)$ we define a leaf $\mathcal{F}_C$  of the absolute period foliation $\mathcal{F}$ and the action of $T_C\EuScript{F}_C\times\mathbb{R}$ on it. This action is described using a geodesic flow associated with the leaves of the foliation $\widehat{\mathcal{F}}$ defined on the stratum of the Teichm\"uller space in the same way as an absolute period foliation. This action is known under the name Rel flow even if it is not formally a flow.

Despite of the existence of numerous results about the absolute period foliations in the case of abelian differentials, almost nothing is known on the analogues of the REL foliations for meromorphic foliations. Recently we were updated about the partial progress obtained by E. Arzhakova, G. Calsamiglia and B. Deroin (\cite{CaDe}, \cite{ArCaDe})
who work on the analogue of the transfer principle applied to meromorphic case with several restrictions on the number of poles or the genus of the surface. It would be interesting to compare their observations with our resutls.

\section{A combinatorial model for the space of generic real-normalized differential}\label{sec:comb}

Let~$\phi$ be a meromorphic differential on a Riemann surface~$C$,
and let~$q\in C$ be a point that is neither a pole, nor a zero of~$\phi$. To these data,
one associates a real line in the tangent line~$T_qC$ to~$C$ at~$q$ in the following way:
this line is tangent to the (germ of the) curve
given by the equation $\Re\int_q^z\phi=0$.
The line is oriented in the direction in which the imaginary part $\Im\int_q^z\phi$ of the integral increases.
The oriented lines defined in this way form a smooth oriented line field on the
complement to the set of poles and zeroes of~$\phi$.
This oriented line field can be extended to the zeroes and poles of~$\phi$ as a singular oriented line field.
We will not use the extension to poles, while
at a zero~$q_0$ of order~$k$ there are $2(k+1)$ distinguished oriented real half-lines in $T_{q_0}C$, whose orientations
alternate around the point~$q_0$.
In particular, at a zero of order~$1$ (a simple zero), the singularity
of the oriented line field looks like two transversally intersecting lines,
the four half-lines being oriented in an alternating order.
We denote this singular oriented line field on the complement to the poles of~$\phi$ by~$I_\phi$.

Now, suppose~$\phi$ is a generic real-normalized meromorphic differential
on~$C$ with a unique pole, of order~$2$, at a point $p\in C$.
Here the term {\it generic\/} means that
\begin{itemize}
\item all the zeroes of~$\phi$
are simple, and
\item the integrals of~$\phi$ between any pair of its zeroes
are not purely imaginary.
\end{itemize}
Note that since all the zeroes are simple, there are exactly~$2g$ of them.

The two integral curves of the oriented line field~$I_\phi$ passing through
each zero of the $1$-form~$\phi$ connect this zero to the only pole of~$\phi$.
Cut the surface~$C$ along the two integral half-curves \emph{entering} each zero of~$\phi$;
the result will be called a \emph{cut surface}.
Picking a point $A$ in the cut surface~$C$ and introducing on this cut surface the
complex coordinate $\int_A^x\phi$,
we identify it with the complex line~$\C$ with $4g$ pairwise distinct distinguished points on it
split into~$2g$ disjoint pairs, which is also cut
along the $4g$ vertical half-lines issuing down from the distinguished points.
The vertical lines $\Re(z)=\const$ form the image of the integral curves of the line field~$I_\phi$ on~$C$,
while the horizontal lines $\Im(z)=\const$ are the integral curves of the second line field, $R_\phi$,
formed by the lines tangent to the germs of the curves
$\Im~\int_q^z\phi=0$, for each point $q\in C$ that is neither a pole,
nor a zero of~$\phi$. Each pair of distinguished points
lies on the same horizontal line. The original surface~$C$ can be reconstructed from the
cut complex line~$\C$ by gluing each left (resp., right) side of a cut to the right
(resp., left) side of the pair cut.

\begin{figure}
  \begin{minipage}[h]{0.5\linewidth}
    \center{\includegraphics[width=0.9\linewidth]{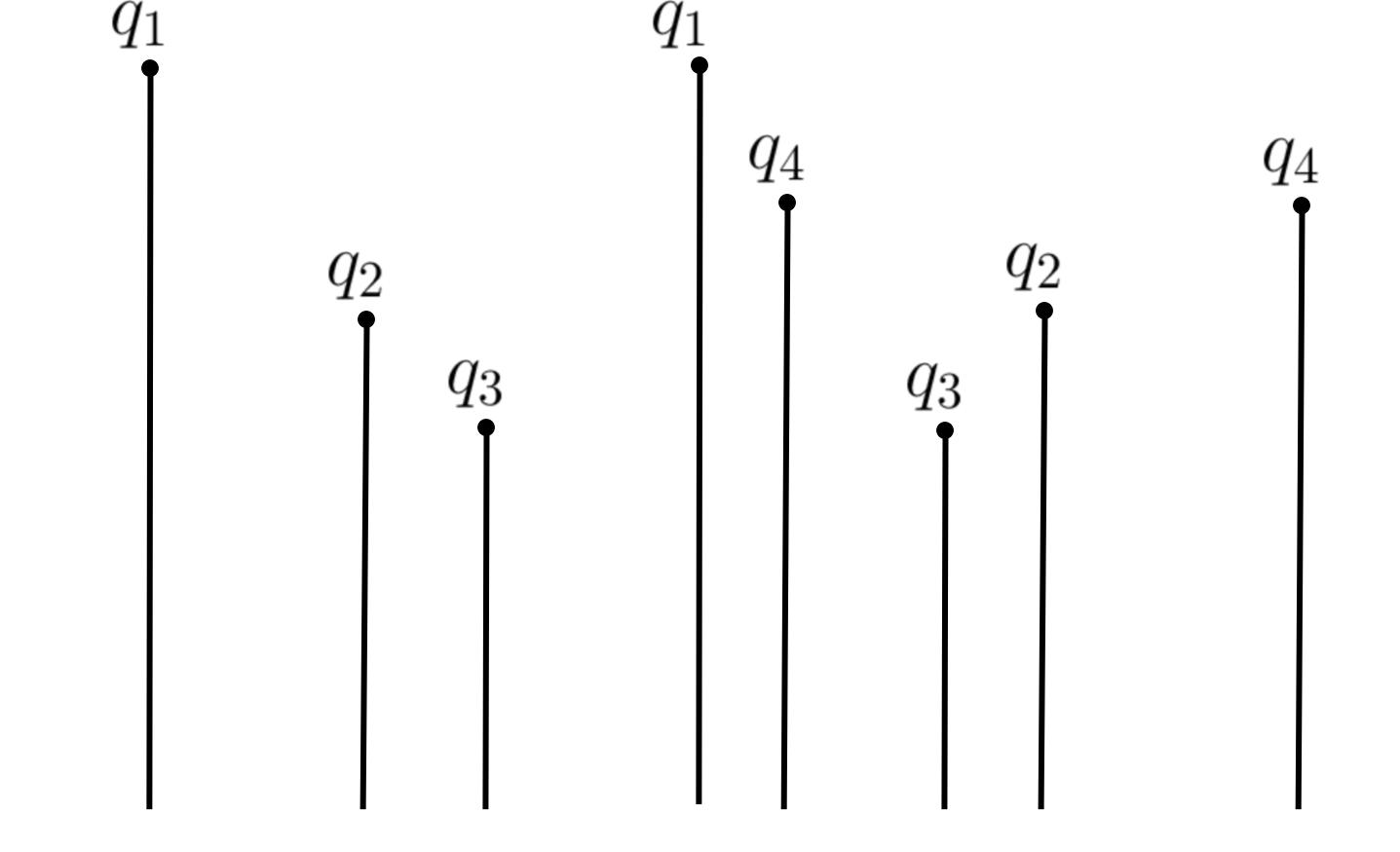}}
  \end{minipage}
  \begin{minipage}[h]{0.5\linewidth}
    \center{\includegraphics[width=0.9\linewidth]{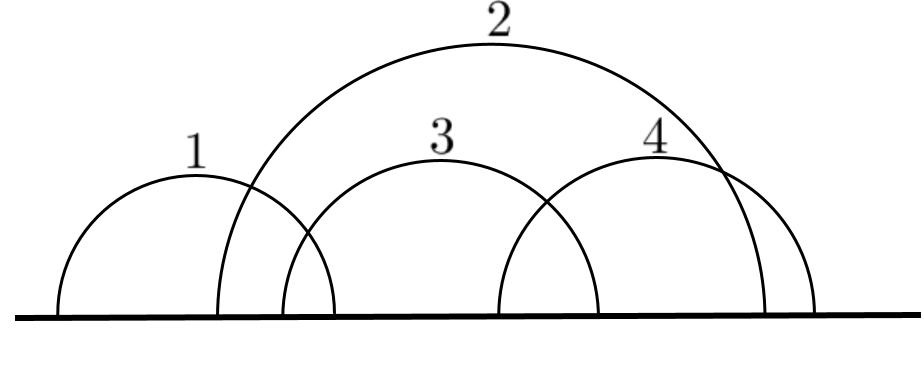}}
  \end{minipage}

  \begin{minipage}[h]{0.5\linewidth}
    \center{(\text{a})}
  \end{minipage}
  \begin{minipage}[h]{0.5\linewidth}
    \center{(\text{b)}}
  \end{minipage}

  \caption{(a) A cut diagram and (b) the corresponding arc diagram
  }\label{fcd}
\end{figure}

In the opposite direction, pick on the complex line~$\C$ any $2g$~pairs $(p_i,q_i)$, $i=1,2,\dots,2g$, of distinct points
 possessing the following properties:
\begin{itemize}
\item in each pair two points have the same imaginary part (but distinct real
parts):
$\Re p_i<\Re q_i$, $\Im p_i=\Im q_i$;
\item points belonging to different pairs have distinct imaginary parts.
\end{itemize}
Now make cuts along vertical half-lines issuing down from all the points $p_i,q_i$
and going to infinity (see Fig.~\ref{fcd}).
To obtain the surface~$C=C_{(p,q)}$, glue the left (respectively, right) side of the cut at the point~$p_i$
to the right (respectively, left) side of the cut at the pair point~$q_i$, identifying points whose coordinates have the same imaginary part.
Under this gluing, the meromorphic $1$-form $dz$ induces a meromoprphic $1$-form on~$C_{(p,q)}$,
which we denote by~$\phi_{(p,q)}$. The $1$-form $\phi_{(p,q)}$ is real-normalized. Indeed,
the group of its periods is generated by the distances $q_i-p_i$, which are real.

\begin{remark}
The shift of all the points $p_i,q_i$, which determine the cutting, by the same complex vector~$\tau$
does change neither the resulting surface $C_{(p,q)}$, nor the $1$-form $\phi_{(p,q)}$ on it.
This shift corresponds to the choice of the base point~$A$ on the cut surface~$C$.
Therefore, we may choose an appropriate normalization
of the set $(p,q)$. One of the possible normalizations can require that
$$
\sum_{i=1}^{2g}(p_i+q_i)=0.
$$
\end{remark}

\begin{remark}
 At a simple zero~$p_0$ of the
	$1$-form $\phi$, there are two integral curves of the line field~$R_\phi$, which is complimentary to~$I_\phi$,
	leaving the zero. In one of the two directions, the integral curve of~$I_\phi$
	approaches infinity, while in the other one it returns to~$p_0$; the integral of~$\phi$
	over this closed curve is a period of~$\phi$. This period is strictly positive.
\end{remark}

If we want that this resulting $1$-form has a unique pole on~$C$, we must impose
a non-degeneracy condition. This non-degeneracy condition is best expressed in terms of the arc
diagram associated to the given cut diagram. Pick a horizontal line in~$\C$
on the level that is so low that each of the $4g$ vertical cuts intersects it.
Then the cuts intersect the horizontal line at $4g$ pairwise distinct points split into
$2g$~pairs. Requiring that the $1$-form $\phi$ has a single pole means that
after cutting the horizontal line and regluing it according to the
gluing scheme given by the cuts, we obtain a connected line once again (and not a disjoint union of a line and a number
of circles).

We connect points belonging to the same pair by an arc in the upper
half-plane and get an object called an  \emph{arc diagram}. We say that two arcs \emph{intersect one another} if their ends
follow the horizontal line in alternating order. The \emph{intersection graph} of an arc
diagram is the graph whose vertices correspond one-to-one to the arcs of the diagram, and
two vertices are connected by an edge iff the corresponding arcs intersect one another.
The following statement justifies the degeneracy condition.

\begin{theorem}\cite{B77,S01}\label{th:ES}
The result of cutting and gluing the horizontal line is a connected line iff the adjacency
matrix of the intersection graph of the arc diagram is nondegenerate over ${\mathbb F}_2$.
\end{theorem}

For example, the adjacency matrix of the intersection graph of the arc diagram in Fig.~\ref{fcd},
which is
$$
\left(\begin{array}{cccc}
0&1&1&0\\
1&0&0&1\\
1&0&0&1\\
0&1&1&0
\end{array}\right),
$$
is degenerate. This means that this arc diagram cannot appear as the result of cutting a
(genus~$1$) surface endowed with a real-normalized differential
having a single pole of order~$2$.

(In fact, the statement proved in~\cite{B77,S01} is more general. It asserts that
{\it the number of connected components of the result of regluing the horizontal line
is one more than the corank of the adjacency matrix of the intersection graph}. This statement describes the number of poles of the meromorphic
differential form obtained from $dz$ on the Riemann surface, which
results from the gluing of~$\C$ along the cuts.
All the resulting poles but one are of order~$1$;
they correspond to compact connected components
of the reglued horizontal line.
The last pole, which corresponds to the segment of
the reglued horizontal line containing both infinite rays,
is of order~$2$.)

Now, the set of nondegenerate arc diagrams with~$2g$ arcs
is in a one-to-one correspondence with the set of connected
components in the space
of generic real-normalized differentials:
any two cut diagrams whose arc diagrams coincide with a given
one can be connected by
a smooth path in the set of cut diagrams having the same arc diagram,
while cut diagrams with different arc diagrams belong to different connected components.
Therefore, counting
non-degenerate arc diagrams with $2g$ arcs allows one
to count connected components in the space of generic real-normalized
differentials.

\begin{proposition}
	The number of connected components in the space of generic real-normalized
	differentials with a single pole of order~$2$ and no other poles on
	genus~$g$ surfaces is equal to
	$$
	\frac{(4g-1)!!}{2g+1}.
	$$
\end{proposition}

 Here the double factorial $(2m-1)!!$ of a positive odd number
 denotes the product of all odd numbers from~$1$ to~$(2m-1)$,
 $(2m-1)!!=1\cdot3\cdot5\cdot\dots\cdot(2m-1)$.
 The sequence of these numbers starts with
 $$
 1,21,1485,225225,\dots
 $$
and coincides with the sequence $\lambda_g(2g)$ in~\cite{HZ},
 which enumerates strata of maximal dimension $4g-3$ in the Strebel stratification of the
 combinatorial moduli space $\cM_{g;1}^{{\rm comb}}=\cM_{g;1}\times\R_+$.

\section{Absolute period foliation of the moduli spaces of real-normalized differentials}\label{sec:princ}

In this section, we introduce the absolute period foliation on the space of real-normalized
differentials with a singke pole of order~$2$
and prove that the given period lattice locus 
of this foliation is dense for almost each period lattice.

For a real-normalized meromorphic differential~$\phi$ on a Riemann surface~$C$,
denote by~$L_\phi$ the additive subgroup in~$\R$ consisting of the periods of~$\phi$
(the integrals of the $1$-form~$\phi$ over closed curves in~$C$ non passing through
its poles). For a generic differential~$\phi$, the subgroup $L_\phi$ coincides with the one consisting of $\Z$-linear
combinations of the numbers $p_i-q_i$, $i=1,2,\dots,2g$.
Indeed, the homology classes of the closed curves that are represented
by arcs in the complement to the cuts in~$\C$ connecting the
points~$p_i$ to~$q_i$ form a basis in $H_1(C;\Z)$,
and the integral of~$dz$ over such an arc is independent of the choice
of the arc, whence $p_i-q_i$.
Note that each subgroup that can be generated by $2g$ (or fewer) reals can appear in this way.
Remark also that there are two different types of these subgroups: either such a subgroup is discrete, generated by a single
element $a\in\R$, $a>0$, or it is dense in~$\R$. The subgroup $L_\phi\subset\R$
is dense iff~$\phi$ has at least two incommensurable periods.
In addition, we say that a real-normalized differential~$\phi$
on a genus~$g$ surface~$C$ is \emph{totally incommensurable} if its
period subgroup $L_\phi\subset\R$ cannot be generated by less than~$2g$ reals, that is, if $L_\phi\subset\R$ is a free $\Z$-module
of rank~$2g$.

Denote by~$\cR_g$ the space of real-normalized meromorphic differentials with a single pole of order~$2$
on genus~$g$ surfaces.
For each additive subgroup $L\subset\R$ generated by at most $2g$ elements, consider
the locus $A_L\subset\cR_{g}$ of the period foliation in~$\cR_{g}$,
which consists of the real-normalized differentials~$\phi$ such that $L_\phi=L$.

The following statement is easy to prove.

\begin{proposition}\label{pr-d}
For a dense period subgroup~$L\subset\R$, the locus $A_{L}\subset\cR_{g}$
is dense in~$\cR_g$.
\end{proposition}

Indeed, if a subgroup $L$ with at most $2g$ generators as a $\Z$-module is dense in~$\R$, then the
subset in~$\R^{2g}$ consisting of vectors
$(a_1,\dots,a_{2g})$ such that the
$2g$-tuple $\{a_1,\dots,a_{2g}\}$ of their coordinates generate~$L$ is dense in~$\R^{2g}$.
Now, if we take a generic~$\phi\in~\cR_g$ and its cut representation, then in any neighborhood of this
cut representation we can find a cut representation of a generic $1$-form having the same
arc diagram as~$\phi$, with the group of periods coinciding with~$L$:
it suffices to replace the $2g$-tuple of absolute periods of~$\phi$ by a sufficiently
close generic $2g$-tuple of numbers generating the group~$L$.

\begin{theorem}\label{th-c}
For a {totally incommensurable} $1$-form~$\phi$, 
the locus $A_{L_\phi}\subset\cR_{g}$
contains at least as many connected components as the set of
conjugacy classes ${\rm GL}(2g,\Z)/{\rm Sp}(2g,\Z)$. In particular, this locus consists of two
connected components for $g=1$, and consists of infinitely
many connected components for $g>1$.
\end{theorem}

Thus, in contrast to the holomorphic case, where a typical isoperiodic
locus is connected, a typical isoperiodic locus in the space of
real-normalized differentials consists of infinitely many connected leaves.


{\bf Proof of Theorem~\ref{th-c}.}
For a given cut diagram, when moving a pair of corresponding cuts vertically (that is, changing the imaginary parts of the cutting points)
we do not change the subgroup $L_\phi$ of periods of the corresponding $1$-form $\phi$.
Therefore, each cut diagram can be connected by a continuous path inside a locus $A_{L_\phi}$
with a cut diagram all whose cut points have the same imaginary part.
Without loss of generality, this imaginary part can be set equal to~$0$.
Similarly, we can move any pair of corresponding cuts horizontally, inside a given locus $A_{L_\phi}$,
while each cut is disjoint from any other cut.

\begin{figure}

  \center{\includegraphics[width=0.75\linewidth]{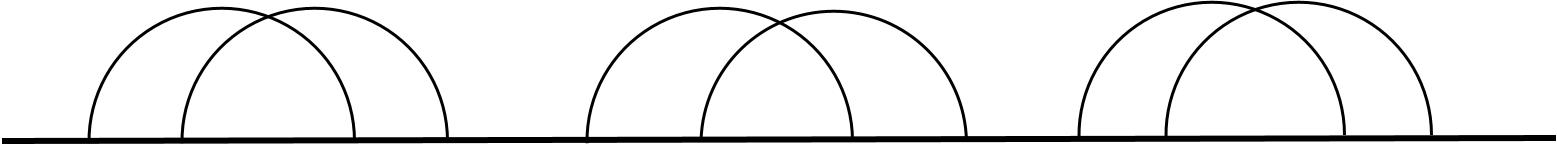}}
\caption{The $3$-caravan
}\label{f-3c}
\end{figure}

By a {\it $g$-caravan}, we will mean the arc diagram formed by~$g$ consecutive pairs of arcs, from left to right,
such that the arcs in each pair intersect one another, while the arcs in different pairs
do not intersect, see Fig.~\ref{f-3c}. A cut diagram whose corresponding arc diagram is a $g$-caravan
will be called a {\it $g$-caravan cut diagram},
and the corresponding cut-surface is a {\it $g$-caravan cut surface}.

Now we are going to prove

\begin{proposition}
For a {totally incommensurable} period subgroup~$L_\phi$, each generic cut surface can be connected by a continuous path in
the period locus  $A_{L_\phi}$
it belongs to to a $g$-caravan cut surface.
\end{proposition}

\begin{figure}
  \center{\includegraphics[width=0.95\linewidth]{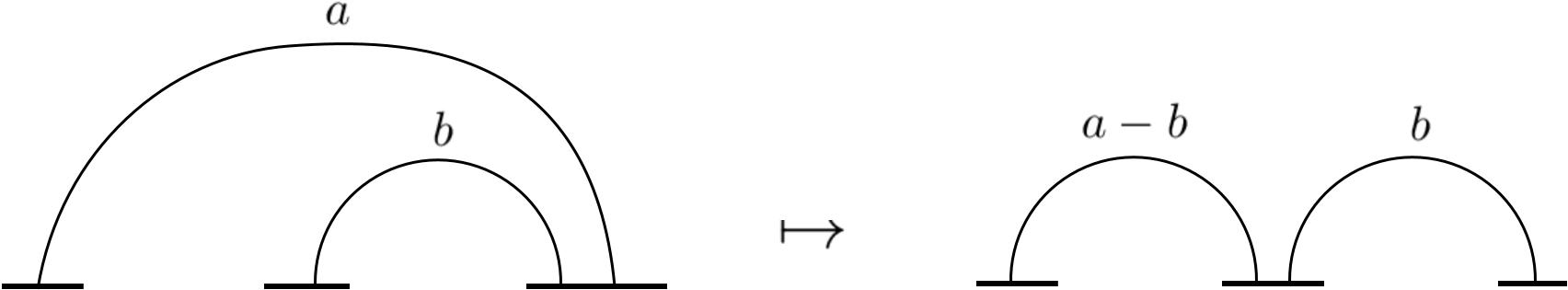}}
\caption{A switch of a cut: as an outer arc of length~$a$ moves to the left, its right end meets
the right end of the inner arc, of length~$b$, $b<a$, and switches to the left end of the inner arc;
as a result the length of the first arc becomes~$a-b$
}\label{f-2Vm}
\end{figure}

Consider a cut diagram, with the imaginary part of all cut points being equal to~$0$.
Moving a pair of corresponding cuts horizontally, say, in the positive direction,
we remain in the same leaf of the period foliation,
and this assertion remains true even if the two cuts coincide (meaning, in particular,
the cut surface fails to be generic). The requirement that the
1-form~$\phi$
is {totally incommensurable}, guarantees that only one pair of cuts, not two pairs, can meet simultaneously
under a horizontal move of a pair of cuts.
As the pair of corresponding cuts
moves further in the same direction, the cut that has met another one,
switches to the pair of this other cut and continues moving, see Fig.~\ref{f-2Vm}.

With this move, the arc diagram of the cut diagram changes undergoing the so called
{\it second Vassiliev move}, which consists in sliding one end of an arc along the other
arc. In~\cite{BNG96}, it is shown that for any nondegenerate arc diagram
of $2g$ arcs there is a sequence of second Vassiliev moves making it into the $g$-caravan.
In our case, the proof must be modified, because of the presence of a metric ---
distances between the arcs' ends, and we present a modification of the proof below.

{\bf Step 1.} Take the arc with the leftmost end and start moving it to the left,
while preserving the positions of the arc ends belonging to all the other arcs, until the right end
of the arc meets an end of another arc. If the two arcs do not intersect one another
(that is, the right end of the first arc meets the right end of the second one),
then apply the second Vassiliev move for the first arc with respect to the second one
(as in Fig.~\ref{f-2Vm}), and repeat the process. Since the arc diagram is nondegenerate,
finally we arrive at the situation where the second arc intersects the first one.
Starting from this point, we will move the two intersecting arcs to the left as a whole
aiming at isolating this pair of arcs.

\begin{figure}
  \center{\includegraphics[width=1\linewidth]{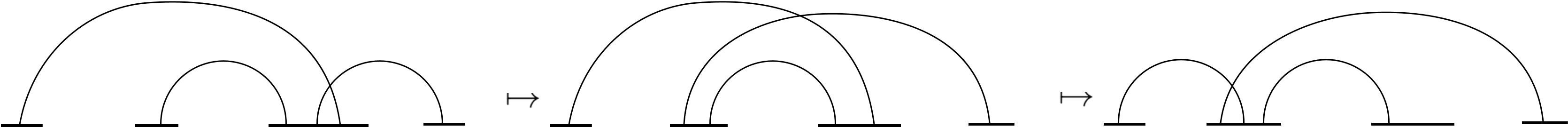}\\ (a)}

   \center{\includegraphics[width=0.66\linewidth]{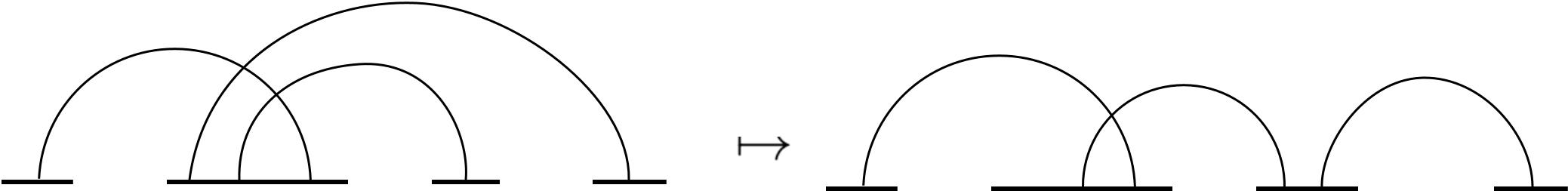}\\ (b)}


    \center{\includegraphics[width=0.66\linewidth]{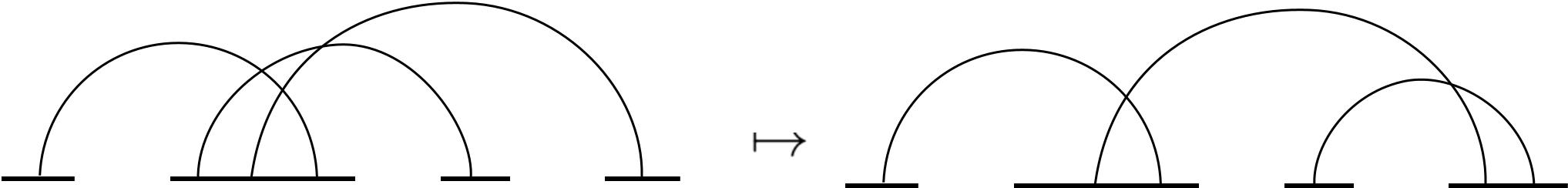}\\ (c)}
\caption{Moving a pair of intersecting arcs, one of which has the leftmost end, to the left,
for the case where the first arc end met belongs to the left interval:
(a) the second end of the arc met also belongs to the left interval; (b) the second end of the
arc met does not belong to the left interval and the arc does not intersect the second one;
(c) the second end of the
arc met does not belong to the left interval and the arc intersects the second one }\label{f-l1m}
\end{figure}

{\bf Step 2.} Now we have two intersecting arcs, one of them with the leftmost end,
and such that the left end of the second arc is the left neighbor of the right end of
the first arc. We will refer to the interval between the left ends of the first and the second arcs as the
left interval, and to the interval between their right ends as the right interval.
As we move two intersecting arcs with neighboring ends to the left,
an end of another arc can meet either the left or the right end of the second arc.

2.1. Consider first the case where the arc end met is in the left interval.
In this case we apply the second Vassiliev move
to the third arc by sliding its end along the second arc, see Fig.~\ref{f-l1m}. As a result, the
number of arc ends in the left interval decreases by one, while the
number of arc ends inside the right interval remains unchanged.

\begin{figure}

  \center{\includegraphics[width=0.75\linewidth]{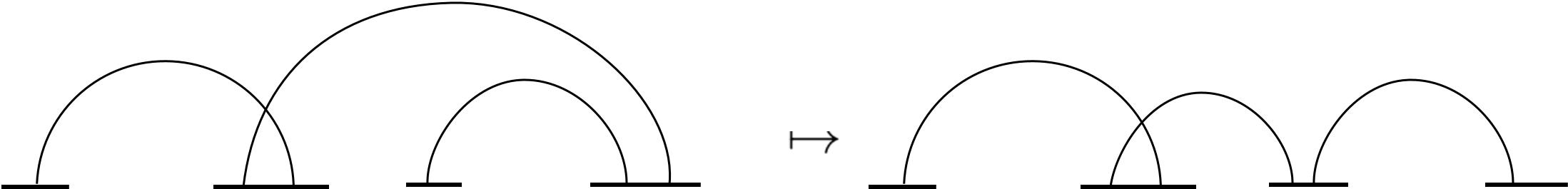}\\ (a)}
  \center{\includegraphics[width=1\linewidth]{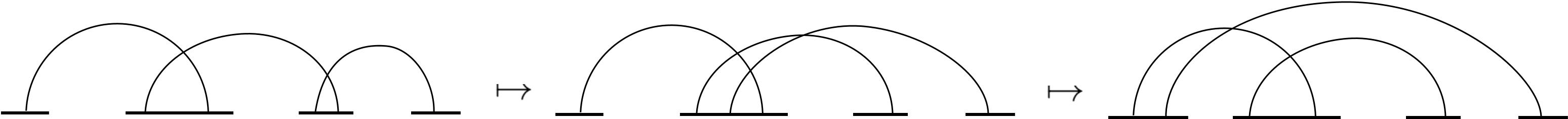}\\ (b)}
\caption{Moving a pair of intersecting arcs, one of which has the leftmost end, to the left,
for the case where the first arc end met belongs to the right interval:
(a) the second end of the arc met also belongs to the right interval; (b) the
arc met intersects the second one}\label{f-l2m}
\end{figure}

2.2. Now suppose the arc end met is in the right interval. If this is
the right end of the third interval, then we apply the second Vassiliev move
by sliding the right end of the second arc along the third arc.
As a result, the
number of arc ends in the right interval decreases at least by two, while the
number of arc ends inside the left interval remains unchanged.
If, in c this is the left end of the third interval (meaning that the second and the
third arcs intersect one another), then by sliding the left end of the third
arc first along the second arc, next along the first arc, we move this left
end from the right interval to the left one. As a result, the number
of arc ends in the right interval decreases by one,
while the number of arc ends in the left interval increases by one, see Fig.~\ref{f-l2m}.

Thus, in each situation, we either decrease the total number of arc ends in both
left and right intervals, or decrease the number of arc ends in the right interval.
Hence, proceeding by induction we can isolate the configuration of the first
and the second arcs.

Now, isolating~$g$ pairs of pairwise intersecting arcs one-by-one,
we obtain a $g$-caravan cut surface connected by a continuous path to
the cut surface we started with, inside the leaf it belongs to.
The Proposition is proved.

To complete the proof of Theorem~\ref{th-c},
remark first that, for a totally incommensurable $1$-form~$\phi$,
any $g$-caravan in $A_{L_\phi}$ defines a basis in the $\Z$-module
$L_\phi$, which consists of the lengths of the arcs, numbered
from left to right. Any two such bases are related by a unique linear
transform in~$GL(2g,\Z)$. On the other hand, any $g$-caravan
determines a basis in the first homology of the underlying $2$-dimensional
surface with coefficients in~$\Z$, and the transition matrix in
this homology from one basis to another one belongs to the group
$\Sp(2g,\Z)$. Hence, the set of connected components of the locus
$A_{L_\phi}$ is at least as large as the set $GL(2g,\Z)/\Sp(2g,\Z)$.
Finally, for $g=1$, any element in~$\Sp(2,\Z)$ can be easily realized
by a sequence of second Vassiliev moves.

\section{Stratification of the moduli spaces of real-nor\-ma\-li\-zed differentials: a case study}\label{sec:str1}

In this section we study the real stratification of the moduli space of real-normalized differentials
with a single pole of order~$2$ on {\it stable} elliptic curves (possibly singular), the pole being a marked point on the curve. 

Recall that meromorphic differentials with a single order $2$ pole
form a rank~$2$ vector bundle over $\overline{\mathcal M}_{g,1}$,
with the fiber over a
stable curve $(C,p)$ being $\omega_C+2p$, where $\omega$ is the relative
dualizing sheaf over $\overline{\mathcal M}_{g,1}$. 
The fiber of this bundle over the singular genus~$g$ curve 
consists of meromorphic differentials over the normalization of
this curve, which is a smooth curve, having
poles of order at most~$1$ at the two preimages of each double point 
(\emph{node}),
with mutually canceling residues.  
In general, taking an analytic section of this bundle over ${\mathcal M}_{g,1}$, one expects that its limit over a singular curve has
simple poles at the nodes. It turns out that this does not happen for real-normalized differentials of the second kind (i.e. having vanishing residues at all poles), and, in particular, for  real-normalized differentials with a single pole of order~$2$.
Recall that the space of real-normalized differentials with a single pole of order~$2$ over genus~$g$ curves is naturally identified with the total
space of the complex line bundle $\cL_1^\vee$ over $\overline{\mathcal M}_{g;1}$.

\begin{theorem}[\cite{GrKr09}]
Any bounded real analytic section  of $\cL_1^\vee$
over~${\mathcal M}_{g;1}$ extends to a {\it continuous}
section of the extension of this bundle to $\overline{\mathcal M}_{g;1}$. 
The limit of the section over a
singular curve~$C$ is the unique meromorphic differential
that is identically zero on all connected components of the normalization of~$C$
except the one containing the pole $p$. 
On that component, it is the unique differential
with real periods and prescribed singular part at $p$.
\end{theorem}
\begin{remark}
Degenerations of real-normalized differential are described
in full generality in \cite{GrKrNo}.
\end{remark}
\bigskip

Now let us return to elliptic curves.
The stratification in question is defined by the mutual
positions of the zeroes of the differential~$\phi$ with respect to the line
field~$I_\phi$. A meromorphic differential on an elliptic curve
with a single pole, the order of which is~$2$, has either two distinct zeroes of order~$1$, or
a single zero of order~$2$. In addition, one of the separatrices of the line field~$I_\phi$,
or both of them, passing through a zero of order~$1$
may contain the other zero of order~$1$.

On a singular elliptic curve, which is the rational curve with two
points glued together, the meromorphic differential~$\phi$ with a single
pole of order~$2$ and no other poles has no zeroes. However,
through each of the two
glued points an integral curve of the line field~$I_\phi$ passes,
and we shall consider the half of this integral curve
entering the point as the separatrix of this point.

From the point of view of the pole,
there are four different cases
depending on the number of separatrices leaving the pole and entering the zeroes:
there could be~$4$ separatrices (the principal stratum), $3$ separatrices
(the stratum in which one of the simple zeroes belongs to a single separatrix entering the other simple zero,
including the case of double zero on a smooth curve), $2$ separatrices
(the stratum in which one of the simple zeroes belongs to both separatrices entering the other simple zero, as well as the stratum of generic real normalized differentials
on the singular elliptic curve),
and $1$ separatrix (the stratum consisting of singular curves
for which a separatrix leaving a preimage
of a double point in the normalization of the singular
curve contains its other preimage).

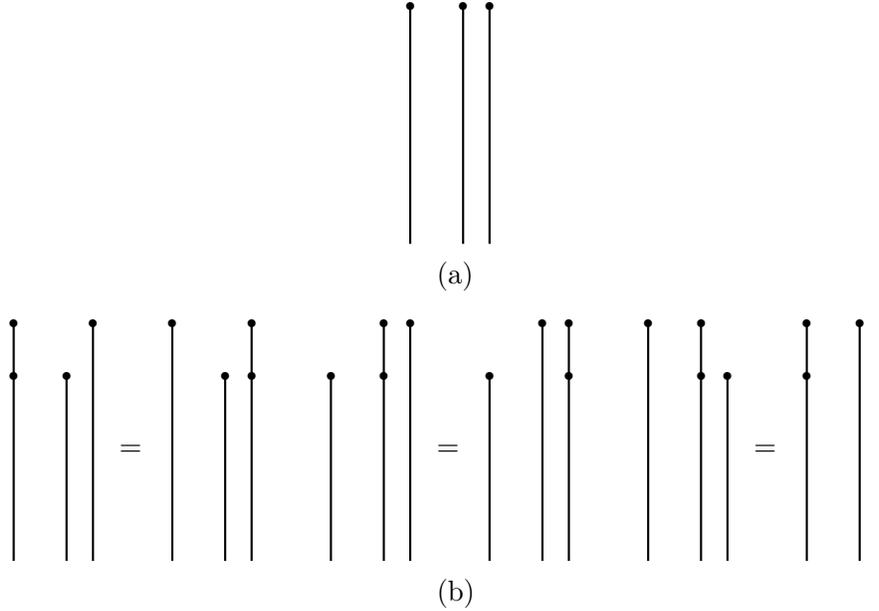
\begin{figure}
\thicklines
\begin{picture}(400,240)(70,-20)
\put(240,220){\line(0,-1){90}}
\put(240,220){\circle*{3}}
\put(260,220){\line(0,-1){90}}
\put(260,220){\circle*{3}}
\put(270,220){\line(0,-1){90}}
\put(270,220){\circle*{3}}
\put(250,115){(a)}
\multiput(90,100)(30,0){2}{\line(0,-1){90}}
\multiput(90,100)(30,0){2}{\circle*{3}}
\put(110,80){\line(0,-1){70}}
\put(90,80){\circle*{3}}
\put(110,80){\circle*{3}}
\multiput(150,100)(30,0){2}{\line(0,-1){90}}
\multiput(150,100)(30,0){2}{\circle*{3}}
\put(170,80){\line(0,-1){70}}
\put(180,80){\circle*{3}}
\put(170,80){\circle*{3}}
\put(130,50){$=$}
\multiput(230,100)(10,0){2}{\line(0,-1){90}}
\multiput(230,100)(10,0){2}{\circle*{3}}
\put(210,80){\line(0,-1){70}}
\put(210,80){\circle*{3}}
\put(230,80){\circle*{3}}
\multiput(290,100)(10,0){2}{\line(0,-1){90}}
\multiput(290,100)(10,0){2}{\circle*{3}}
\put(270,80){\line(0,-1){70}}
\put(300,80){\circle*{3}}
\put(270,80){\circle*{3}}
\put(250,50){$=$}
\put(250,-5){(b)}
\multiput(330,100)(20,0){2}{\line(0,-1){90}}
\multiput(330,100)(20,0){2}{\circle*{3}}
\put(360,80){\line(0,-1){70}}
\put(350,80){\circle*{3}}
\put(360,80){\circle*{3}}
\multiput(390,100)(20,0){2}{\line(0,-1){90}}
\multiput(390,100)(20,0){2}{\circle*{3}}
\put(420,80){\line(0,-1){70}}
\put(420,80){\circle*{3}}
\put(390,80){\circle*{3}}
\put(370,50){$=$}
\end{picture}
\caption{(a) An element in~$\cR^0_1(3)$
	and (b) three pairs of equivalent cut diagrams in $\cR_{1}(3)$ associated to this element
}\label{f-g13}
\end{figure}

\begin{figure}
  \thicklines
  \begin{picture}(400,240)(70,-20)
    \multiput(240,220)(0,-10){10}{\line(0,-1){5}}
    \put(240,220){\circle*{3}}
    \multiput(270,220)(0,-10){10}{\line(0,-1){5}}
    \put(270,220){\circle*{3}}
    \put(250,115){(a)}
    \multiput(120,100)(30,0){2}{\line(0,-1){30}}
    \multiput(120,100)(30,0){2}{\circle*{3}}
    \multiput(120,70)(30,0){2}{\circle*{3}}
    \multiput(120,65)(0,-10){6}{\line(0,-1){5}}
    \multiput(150,65)(0,-10){6}{\line(0,-1){5}}
    \multiput(200,100)(30,0){2}{\line(0,-1){90}}
    \multiput(200,100)(30,0){2}{\circle*{3}}
    \multiput(205,70)(30,0){2}{\circle*{3}}
    \multiput(205,70)(30,0){2}{\line(0,-1){60}}
    \put(170,-5){(b)}
    \multiput(286,110)(8,0){2}{\line(0,-1){105}}
    \multiput(286,110)(8,0){2}{\circle*{3}}
    \multiput(290,90)(30,0){2}{\line(0,-1){85}}
    \multiput(290,90)(30,0){2}{\circle*{3}}
    \multiput(350,110)(0,-10){11}{\line(0,-1){5}}
    \put(350,110){\circle*{3}}
    \multiput(380,90)(0,-10){9}{\line(0,-1){5}}
    \put(380,90){\circle*{3}}
    \multiput(410,110)(30,0){2}{\line(0,-1){105}}
    \multiput(410,110)(30,0){2}{\circle*{3}}
    \multiput(436,90)(8,0){2}{\line(0,-1){85}}
    \multiput(436,90)(8,0){2}{\circle*{3}}
    \put(360,-5){(c)}
  \end{picture}
  \caption{(a) An element in~$\cR^0_1(2)$;
    (b) a cut diagram representing a smooth
    elliptic curve in  $\cR_{1}(2)$ associated to this element (left) and a generic cut diagram close to it (right)	
    (c) (center) a cut diagram representing a singular
    elliptic curve in  $\cR_{1}(2)$ associated to this element (the dashed lines represent the
    separatrices that are not cuts) together with two close
    generic cut diagrams (left and right)
  }\label{f-g12}
\end{figure}
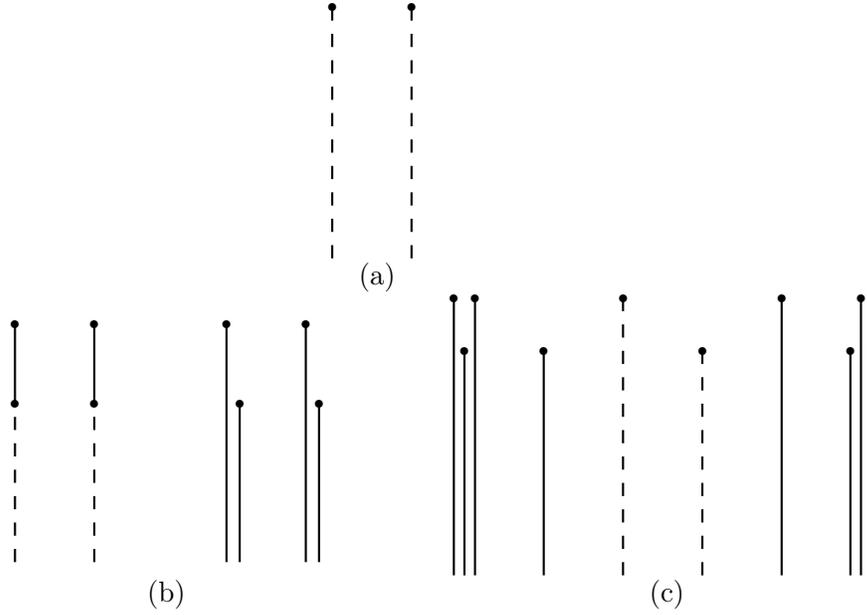

\begin{figure}
  \thicklines
  \begin{picture}(400,140)(70,-20)
    %
    %
    \multiput(200,110)(0,-10){10}{\line(0,-1){5}}
    \put(200,110){\circle*{3}}
    \put(200,90){\circle*{3}}
    %
    %
  \end{picture}
  \caption{A cut diagram representing a singular
    elliptic curve in  $\cR_{1}(1)$
    endowed with a meromorphic differential such that the
    separatrix passing through one of the singular points contains
    the other one
  }\label{f-g11}
\end{figure}
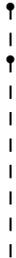

In this way, the moduli space of nonzero real-normalized differentials $\cR_{1}$
is split into the following strata. The cut surfaces corresponding to each of the strata 
are shown in Figures~\ref{f-g14}--\ref{f-g11}.


Here, the identification of the cut surface with the cut diagram is established as follows: take
an arbitrary point~$A$ on the cut surface not belonging to the cuts
and introduce the global coordinate by considering the function $\int_A^q\phi$
of the varying point~$q$ on the cut surface. This coordinate is defined uniquely up
to shifts of the origin determined by the choice of the base point~$A$.
We fix this shift by requiring that the sum of all cut ends is~$0$.

\begin{figure}
	\thicklines
	\begin{picture}(400,240)(70,-20)
	\multiput(220,220)(50,0){2}{\line(0,-1){90}}
	\multiput(220,220)(50,0){2}{\circle*{3}}
	\put(210,222){$p_1$}
		\put(260,222){$q_1$}
	\multiput(250,220)(70,0){2}{\line(0,-1){90}}
\multiput(250,220)(70,0){2}{\circle*{3}}
\put(240,222){$p_2$}
\put(310,222){$q_2$}
	\put(260,115){(a)}
	\multiput(220,90)(50,0){2}{\line(0,-1){90}}
	\multiput(220,90)(50,0){2}{\circle*{3}}
	\put(210,92){$p_1$}
	\put(260,92){$q_1$}
	\multiput(250,70)(70,0){2}{\line(0,-1){70}}
	\multiput(250,70)(70,0){2}{\circle*{3}}
	\put(240,72){$p_2$}
	\put(310,72){$q_2$}
	\put(260,-15){(b)}
	\end{picture}
	\caption{(a) A generic cut surface in~$\cR_1^{0}$ and (b) a
		corresponding element of the
		principal stratum
	}\label{f-g14}
\end{figure}

\begin{proposition}
The general stratum $\cR_{1}(4)$ is a contractible cone $\R_+^3\times\R$ of dimension~$4$ over the $3$-dimensional space $\cR_1^{0}(4)=\R_+^3$.
\end{proposition}

The $\R_+$-cone structure of a general stratum is given by multiplication of a
$1$-form by a positive real constant.
The imaginary part of the coordinates of each copy of the zeroes of~$\phi$
in the cut diagram is the same; we call it the imaginary part of the zero.

To each general cut diagram one can uniquely associate a cut diagram with~$0$
imaginary parts of the zeroes. This mapping defines a fibration of the
space~$\cR_{1}(4)$ of general cut diagrams over the space~$\cR^0_{1}(4)$ of
general cut diagrams with~$0$ imaginary parts of the zeroes.

Now, let $a,b,c$ denote the real parts of the relative periods between the two zeroes,
$a>0,b>0,c>0$. These values can be arbitrary, whence the space~$\cR^0_{1}(4)$
is naturally identified with~$\R_+^3$, so that~$\cR_{1}(4)$ is an~$\R$-fibration
over~$\R_+^3$.

Now, let's turn to the case of three separatrices. The stratum~$\cR_{1}^0(3)$
is the one consisting of real normalized meromorphic differentials having
a single zero, one of order~$2$. Its complement
$\cR_{1}(3)\setminus\cR_{1}^0(3)$ consists of real normalized differentials
with two zeroes, each of order~$1$, one of which belongs to
exactly one of the separatrices entering the second zero.

\begin{proposition}
The stratum~$\cR_{1}^0(3)$ has dimension~$2$
and is homeomorphic to~$\R_+^2$.
The stratum~$\cR_{1}(3)$ is a disjoint union of three copies of $\R_+^2\times\R_+$ of dimension~$3$
and the stratum~$\cR_{1}^0(3)$.
\end{proposition}

Once again, the $\R_+$-cone structure is defined by multiplication of a $1$-form
by a positive real constant.
The space~$\cR^0_{1}$ consisting of cut diagrams with zero imaginary parts
of all the three cut points is the positive real quadrant~$\R_+^2$, the coordinates
$a>0,b>0$ being the two $\phi$-periods. The three $3$-cells $\R_+^3$ attached to
$\cR^0_{1}$, together forming $\cR_{1}$, are shown in Fig.~\ref{f-g13}.

In all the three pairs we choose the left cut diagram as the \emph{canonical} representative of the pair. It is distinguished from the other one by the
requirement that moving from left to right along the low horizontal line,
when meeting a double cut we jump to the nearest pair cut.

The following proposition describes the way the stratum~$\cR_{1}(4)$ is attached to the
stratum~$\cR_{1}(3)$. This description can be given on the level of the strata
$\cR^0_{1}(4)$ and~$\cR^0_{1}(3)$, respectively.

\begin{proposition}
The stratum~$\cR^0_{1}(4)$ is triply attached to the stratum~$\cR^0_{1}(3)$,
depending on around which of the cuts at the double zero contains two simple zeroes, when perturbed.
\end{proposition}

Now, if we have two separatrices, then none of them is really ``separating''
at infinity, meaning that the integral $\int_P^q\phi$ considered as a function
in the varying point~$q$ of the curve is continuous. We shall denote the corresponding
separatrices (or their parts) by dashed lines.

\begin{proposition}
The stratum $\cR_{1}^0(2)$
consists of singular elliptic curves endowed
with a meromorphic differential $dz$ having
a pole of order~$2$ and no other poles,
with two real points with opposite coordinates glued
together. It
is one-dimensional and homeomorphic to the positive half-line~$\R_+$. Its complement in $\cR_1(2)$ consists of two connected
components.
The first one of these connected components
consists of smooth
elliptic curves endowed with a meromorphic differential
having two distinct zeroes, one of which belongs to
both separatrices of the vector field $I_\phi$
entering the second zero and is homeomorphic to
$\R_+\times\R_+$.
The second connected component consists of singular
elliptic curves endowed with a meromorphic differential~$dz$, which
has a single pole of order~$2$ and no other poles,
with two distinct points having opposite coordinates whose real
part is nonzero glued together and homeomorphic to
two copies of $\R_+\times\R_+$.
\end{proposition}

Indeed, if an elliptic curve is singular, then the
integral of the differential form~$\phi$ defines
a coordinate on the rational curve, which is the normalization of the elliptic one. This coordinate
is well defined up to a shift, and this shift can
be chosen so as to make the sum of the coordinates
of the two points glued together equal to~$0$.
If these coordinates are real, then we obtain a point
in~$\cR_{1}^0(2)$. If the original elliptic curve
is smooth and one of the zeroes of the meromorphic
differential belongs to both separatrices of the second
zero, then both the curve and the differential are
totally determined by two positive real numbers,
both being the periods of~$\phi$.

\begin{proposition}
	The stratum $\cR_{1}^0(1)$ is empty. The stratum
	$\cR_1(1)$ consists of singular elliptic curves
	endowed with a meromorphic $1$-form $\phi$
	having a single pole of order~$2$ and no other poles
	such that the coordinates of the two preimages of
	the double point under normalization have coinciding
	real parts. The stratum $\cR_1(1)$ is one-dimensional
	and is homeomorphic to~$\R_+$.
\end{proposition}

\end{document}